\documentclass{article}
\title{The singular locus of Lauricella's $F_C$}
\author{Ryohei Hattori and Nobuki Takayama}
\date{October 29, 2011}
\usepackage{color}
\def\pd#1{ \partial_{#1} }
\def\QED{ Q.E.D.  \bigbreak}
\def\www{{({\bf 0},{\bf 1})}}
\def\km{{(-{\bf 1},{\bf 1})}}   
\def\outertensor{ \mbox{\fbox{$\times$}} }

\newtheorem{theorem}{Theorem}[section]
\newtheorem{proposition}[theorem]{Proposition}
\newtheorem{example}[theorem]{Example}
\newtheorem{corollary}[theorem]{Corollary}

\usepackage{amssymb}
\usepackage{amsmath}
\usepackage{bm}

\DeclareMathOperator*{\ssum}{\bm{\Sigma}}
\DeclareMathOperator*{\pprod}{\bm{\Pi}}

\newcommand{\C}{\mathbb{C}}      
\newcommand{\Z}{\mathbb{Z}}      
\newcommand{\N}{\mathbb{N}}      \newcommand{\di}{\displaystyle}  
\newcommand{\e}{\varepsilon}     

\begin{document}
\maketitle

\noindent
{\bf Abstract.} We determine the singular locus of the Lauricella function 
$F_C$ by utilizing the theory of $D$-modules and Gr\"obner basis. The 
$A$-hypergeometric system associated to $F_C$ is also discussed. 

\section{Introduction}

Let $D=\C\langle x_1, \ldots, x_m, \pd{1}, \ldots, \pd{m} \rangle $ 
be the Weyl algebra of $m$ variables.
We take $2m$ dimensional integer vector $(u,v)$, $u, v \in \Z^m$ 
such that $u_i + v_i > 0$.
For an element 
$p=\sum_{(\alpha, \beta) \in E} c_{\alpha,\beta} x^\alpha \pd{}^\beta$ of $D$, 
we define the $(u,v)$-initial form 
${\rm in}_{(u,v)}(p)$
of  $p$ by
the sum of the terms in $p$ which has the highest $(u,v)$-weight.
In other words, we define
\begin{eqnarray*}
 {\rm ord}_{(u,v)}(p) &=& {\rm max}_{(\alpha,\beta) \in E} 
 (\alpha\cdot u + \beta\cdot v), \\
   {\rm in}_{(u,v)}(p) &= &
 \sum_{(\alpha, \beta) \in E, 
\alpha\cdot u + \beta\cdot v={\rm ord}_{(u,v)}(p)}
  c_{\alpha\beta} x^\alpha \xi^\beta.
\end{eqnarray*}
Here, $\xi_i$ is a new variable which commutes with the other variables
(see, e.g., \cite[\S 1.1]{SST}).
When $u_i+v_i = 0$, we define the $(u,v)$-initial form analogously
and $\xi_i$ is replaced by $\pd{i}$ in the definition above. 
Put ${\bf 0} = (0, \ldots, 0) \in \Z^m$ 
and
${\bf 1}=(1, \ldots, 1) \in \Z^m$.
For $p \in D$, the $({\bf 0},{\bf 1})$-initial form of $p$ is called
the principal symbol of $p$.
For a given left ideal $I$ of $D$,
the characteristic ideal of $I$ is 
${\rm in}_{({\bf 0},{\bf 1})} (I)$,
which is the ideal  in $\C[x_1, \ldots, x_m, \xi_1, \ldots, \xi_m]$
generated by all principal symbols of the elements of $I$.

The zero set in $\C^{2m}$ of the characteristic ideal 
is called the characteristic variety, which is denoted by ${\rm Ch}(I)$.
When the (Krull) dimension of the characteristic variety is $m$, 
then  $D/I$ is called a holonomic $D$-module and
the $I$ is called a holonomic ideal.
The projection of ${\rm Ch}(I) \setminus V(\xi_1, \ldots, \xi_m)$
to the first $m$-coordinates $\C^m = \{ x \}$ 
is called the singular locus of $I$,
which is denoted by ${\rm Sing}(I)$ or
${\rm Sing}(D/I)$.
As to these fundamental notions of $D$-modules, see, e.g., \cite{coutinho}.

Let $R=\C(x_1, \ldots, x_m)\langle \pd{1}, \ldots, \pd{m}\rangle$ 
be the ring of differential operators with rational function coefficients.
The holonomic rank of $I$ is
the dimension of
$R/(RI)$ regarded as 
the $\C(x)=\C(x_1, \ldots, x_m)$ vector space.
The rank is denoted by ${\rm rank}\, (I)$.
The rank is equal to the multiplicity of the characteristic ideal
at a generic point.
In other words, we have
$$ {\rm rank}(I) = {\rm dim}_\C {\cal O}_a
  /{\cal O}_a \cdot {\rm in}_{({\bf 0},{\bf 1})}(I)|_{{x=a}},  
$$
where ${\cal O}_a=\C\{\xi_1-a_1, \ldots, \xi_m-a_m\}$ and 
$a$ is a point in $\C^m \setminus {\rm Sing}(I)$.
We also have the identity
$$
{\rm rank}(I) =
 {\rm dim}_{\C(x)}  \C(x)[\xi]/
  \C(x)[\xi] \cdot {\rm in}_{({\bf 0},{\bf 1})}(I)
$$
where $\C(x)[\xi]$ denotes $\C(x_1, \ldots, x_m)[\xi_1, \ldots, \xi_m]$.
Let ${\rm Sol}(I)$ be the constructive sheaf of
holomorphic solutions on $\C^m$;
$${\rm Sol}(I) = \{ f \in {\cal O} \,|\, \ell \cdot f = 0
  \mbox{ for all } \ell \in I \}.
$$
The holonomic rank ${\rm rank}(I)$ is equal to
${\rm dim}_\C {\rm Sol}(I)(U)$ 
for any simply connected open set $U$ in $\C^m \setminus {\rm Sing}(I)$.
As to these characterizations of the holonomic rank, 
see, e.g., \cite{oaku-char}, \cite[Chapter 1]{SST} and their references.

Put $\theta_i=x_i \pd{i}$ and
$\theta=\sum_{j=1}^m \theta_i$.
We consider the left ideal $I(m)$ generated by the operators
\begin{equation} \label{eq:LauricellaOp}
 \ell_i=\theta_i (\theta_i+c_i-1)  - x_i (\theta+a) (\theta+b), 
\quad i=1, \ldots, m, 
\end{equation}
where $a,b,c_i \in \C$ are parameters.
The Lauricella function $F_C$ is annihilated by the left ideal.

We will show, in Theorem \ref{th:singLocus}, that the singular locus of $I(m)$
agrees with the zero set of
\begin{equation} \label{eq:singLocus}
\prod_{i=1}^m x_i 
  \prod_{\varepsilon_i \in \{ -1, 1 \} } 
    \left( 1 + \varepsilon_1 \sqrt{x_1}+\cdots+\varepsilon_m \sqrt{x_m}\right).
\end{equation}
Note that when we expand (\ref{eq:singLocus}), then it becomes 
a polynomial in $x$.
The proof of this fact occupies the first three sections of this paper.
In the last section, we study the $A$-hypergeometric system
associated to the Lauricella $F_C$
and determine the singular locus of it in the complex torus
by utilizing our main theorem \ref{th:singLocus}.

We had believed 
that the singular locus of $I(m)$ is well-known among experts, 
but we find a few literatures on rigorous proofs on these facts.
Some of them are master thesises by Kaneko \cite{kaneko}
and by Yoshida \cite{yoshida-noumi},
who proved
that the singular locus is contained in (\ref{eq:singLocus})
but they did not discuss on the opposite inclusion.

We thank to Keiji Matsumoto for posing the problem 
and for the encouragement. 

\section{A variety containing the singular locus}

The singular locus of the system $I(m)$ is, by definition,
$\pi \left(V({\rm in}_{({\bf 0},{\bf 1})}(I(m))
\setminus V(\xi_1, \ldots, \xi_m)\right) $.
Here, $\pi$ is the projection from $(x,\xi)$ to $x$.
The principal symbol of $\ell_i$ is denoted by $L_i$ and
they are equal to 
$$ L_i= x_i^2 \xi_i^2 - x_i \left( \sum_{j=1}^m x_j \xi_j \right)^2. $$
In other words, ${\rm in}_{({\bf 0},{\bf 1})} (\ell_i) = L_i$.
Since $L_i \in {\rm in}_{({\bf 0},{\bf 1})}(I(m))$,
the singular locus of the system $I(m)$
is contained in 
$C = \pi (V(L_1, \ldots, L_m) \setminus V(\xi_1, \ldots, \xi_m))$.

Let us regard 
$V(L_1, \ldots, L_m)$ as an analytic space.
When $x_i \not= 0$, $L_i$ is factored as
\begin{equation} \label{eq:chvar}
 L_i = \left( x_i \xi_i - \sqrt{x_i} \left( \sum_{j=1}^m x_j \xi_j \right)
         \right)
         \left( x_i \xi_i + \sqrt{x_i} \left( \sum_{j=1}^m x_j \xi_j \right)
         \right). 
\end{equation}
Therefore, the necessary and sufficient condition that
$x$ lies in $C \cap (\C^*)^m$ is that
\begin{equation} \label{Sakamoto}
 x_i \xi_i +\varepsilon_i \sqrt{x_i} \left( \sum_{j=1}^m x_j \xi_j \right)
=0, \ 
 i=1, \ldots, m
\end{equation}
has a non-trivial solution $\xi\not=0$
for a set of signs $\varepsilon_i$.
The condition can be written in terms of the determinant
of the system regarded as a system of linear equations with respect to
$\xi$.

\begin{proposition}  \label{prop:det}
The determinant of the coefficient matrix of the system of linear equation 
(\ref{Sakamoto}) is equal to
$\di\pprod_{i=1}^mx_i\left( \di 1+\ssum_{j=1}^m\e_j\sqrt{x_j} \right)$.
\end{proposition}

{\it Proof}\/.
The coefficient matrix of the system (\ref{Sakamoto}) is 
\begin{equation}\nonumber
M=\left( 
\begin{array}{ccccc}
\!x_1+\e_1x_1\sqrt{x_1}\!& \e_1x_2\sqrt{x_1}& \e_1x_3\sqrt{x_1}&\cdots&\e_1x_m
\sqrt{x_1}\\
\e_2x_1\sqrt{x_2}&\!x_2+\e_2x_2\sqrt{x_2}\!&\e_2x_3\sqrt{x_2}&\cdots&\e_2x_m
\sqrt{x_2}\\
\e_3x_1\sqrt{x_3}& \e_3x_2\sqrt{x_3}&\!x_3+\e_3x_3\sqrt{x_3}\!&\cdots&\e_3x_m
\sqrt{x_3}\\
\vdots&\vdots&\vdots&\ddots&\vdots\\
\e_mx_1\sqrt{x_m}&\e_mx_2\sqrt{x_m}&\e_mx_3\sqrt{x_m}&\cdots&\!x_m+\e_mx_m
\sqrt{x_m}\!
\end{array}\right). 
\end{equation}
We have 
\[\begin{array}{cl}
 &\det M\\
=& x_1x_2\cdots x_n\det
\left( 
\begin{array}{ccccc}
\!1+\e_1\sqrt{x_1}\!& \e_1\sqrt{x_1}& \e_1\sqrt{x_1}&\cdots&\e_1
\sqrt{x_1}\\
\e_2\sqrt{x_2}&\! 1+\e_2\sqrt{x_2}\!&\e_2\sqrt{x_2}&\cdots&\e_2
\sqrt{x_2}\\
\e_3\sqrt{x_3}& \e_3\sqrt{x_3}&\!1+\e_3\sqrt{x_3}\!&\cdots&\e_3
\sqrt{x_3}\\
\vdots&\vdots&\vdots&\ddots&\vdots\\
\e_m\sqrt{x_m}&\e_m\sqrt{x_m}&\e_m\sqrt{x_m}&\cdots&\!1+\e_m\sqrt{x_m}\!
\end{array}\right)\\
=& \di\pprod_{i=1}^m x_i\det\left( 
\begin{array}{crrrr}
1+\e_1\sqrt{x_1}& -1& -1&\cdots&-1\\
\e_2\sqrt{x_2}& 1&0&\cdots&0\\
\e_3\sqrt{x_3}&0&1&\cdots&0\\
\vdots&\vdots&\vdots&\ddots&\vdots\\
\e_m\sqrt{x_m}&0&0&\cdots&1
\end{array}\right)
\left( \begin{array}{l}
{[j \mathrm{th~column}]}\\
\!\!{+[1\mathrm{st~column}]\times(-1)}\!\! \\
{\hspace{4.2pt}\mathrm{for~} 2\leq j\leq m}
\end{array}\right)
\\
=& \di\pprod_{i=1}^m x_i\det\left( \begin{array}{cc}
1+\di\ssum_{j=1}^m\e_j\sqrt{x_j}&0\\
*& E_{m-1}\end{array}\right)\quad \left( {[1\mathrm{st~row}]+\di\ssum_{j=2}
^m[j\mathrm{th~row}]}
\right)\\
=& \di\pprod_{i=1}^mx_i\left( \di 1+\ssum_{j=1}^m\e_j\sqrt{x_j}\right)
\end{array}\]
\QED

\begin{theorem} {\rm \cite{kaneko}, \cite{yoshida-noumi}}
The singular locus of $I(m)$ is contained in the zero set of
(\ref{eq:singLocus}).
\end{theorem}

{\it Proof}\/.
Since $x_i=0$ are already the zero set of (\ref{eq:singLocus}), we may only 
consider the singular locus
in $(\C^*)^m$.
If $x \in (\C^*)^m \cap C$, then
(\ref{Sakamoto}) must have a non-trivial solution $\xi\not=0$.
By Proposition \ref{prop:det}, we get the conclusion.
\QED

In the sequel, we want to prove the opposite inclusion
$ C \subseteq {\rm Sing}(I)$
where $C$ is the zero set of (\ref{eq:singLocus}).
If a classical solution of $I(m)$ has singularities on the all
irreducible components of the zero set $C$,
then we are done.
However, as the following example shows,
the singular locus of classical solutions may smaller than
the zero set $C$.

\begin{example} \rm
Assume $m=2$.
When $a=-1/2, b=-2, c_1=c_2=1/2$,
the solution space of the differential equation
is spanned by the following functions
$$
1+2 x+2 y-2 x y-x^2/3-y^2/3,
\sqrt{x}, \sqrt{y},
\sqrt{xy} (1-x/3-y/3)
$$
Note that the singular locus of these solutions is contained in
$xy=0$,
which is smaller than the zero set $C$.
\end{example}

We close this section with two preparatory propositions.
\begin{proposition}
The left ideal $I(m)$ is holonomic.
\end{proposition}

{\it Proof}\/.
Since the Bernstein inequality ${\rm dim}\, V({\rm in}_\www(I(m))) \geq m$ 
holds (see, e.g., \cite{coutinho}, \cite{SST}),
we may prove that ${\rm dim}\, V({\rm in}_\www(I(m))) \leq m$.
It follows from the decomposition (\ref{eq:chvar})
that for any $x \in (\C^*)^m$
the $\xi$'s satisfy (\ref{eq:chvar}) are finite points. 
Hence, the analytic set $C \cap (\C^*)^m$
is $m$-dimensional.
The remaining thing to do is the evaluation of the dimension at the points
in $x_i=0$.

We put $I'(m) = \langle L_1, \ldots, L_m \rangle$.
$I'(m)$ is contained in ${\rm in}_\www(I(m))$.
We will prove 
${\rm dim}\, I'(m) = m$ by induction.
When $m=1$, it is easy to see that ${\rm  dim}\, I'(m) = 1$.
Let us assume ${\rm dim}\, I'(m-1) = m-1$.
We note that $V(I'(m)) \cap V(x_{m}) = 
\{ ((x',\xi'), (0,\xi_m)) \,|\, (x',\xi') \in V(I'(m-1)) \subset \C^{m-1},
  \xi_m \in \C \}$
because $x_m \xi_m = 0$ in $L_i$ when $x_m = 0$. 
It follows from the induction hypothesis 
${\rm dim}\, V(I'(m-1)) = m-1$ that
the dimension of $V(I'(m))$ at any point in $x_m=0$ is $(m-1)+1=m$.
\QED

\begin{proposition}
The polynomial 
$$  \prod_{\varepsilon_i \in \{+1,-1\}}   \left( 1 + \varepsilon_1 \sqrt{x_1}+\cdots+\varepsilon_m \sqrt{x_m}\right)$$
is irreducible
in $\C[x_1, \ldots, x_m]$.
\end{proposition}

{\it Proof}\/.
Let 
\begin{eqnarray*}
P(t) &=& P(t_1,\ldots,t_m) ~=~ \prod_{\varepsilon\in \{ \pm 1\}^m}P_\varepsilon
(t),\\
P_\varepsilon(t) &=& 1+\e_1t_1+\e_2t_2+\cdots+\e_mt_m, 
\end{eqnarray*}
where $\e=(\e_1,\ldots,\e_m)$. 
First we show that the polynomial $P(t)$ is irreducible in the subring 
$\C[t_1^2,t_2^2,\ldots,t_m^2]$ of $\C[t]$. 
Suppose that there exist $Q_1(t), Q_2(t)\in \C[t]$ satisfying 
\[ P(t)~=~ Q_1(t_1^2,t_2^2,\ldots,t_m^2)Q_2(t_1^2,t_2^2,\ldots,t_m^2).\]
Then in $\C[t]$, we have in particular 
\[ P_{\bm{1}}(t)\mid 
Q_1(t_1^2,t_2^2,\ldots,t_m^2)Q_2(t_1^2,t_2^2,\ldots,t_m^2),\]
where $\bm{1}=(1,1,\ldots,1)\in \{ \pm 1\}^m$. 
Since $P_{\bm{1}}(t)=1+t_1+t_2+\cdots+t_m$ is a prime element in $\C[t]$, 
either $Q_1(t_1^2,t_2^2,\ldots,t_m^2)$ or $Q_2(t_1^2,t_2^2,\ldots,t_m^2)$ is 
divided by $P_{\bm{1}}(t)$ in $\C[t]$. Assume the former, then we have 
\begin{equation}\label{P1divQ1}
Q_1(t_1^2,t_2^2,\ldots,t_m^2) = P_{\bm{1}}(t)g(t)
\end{equation}
for some $g(t)\in \C[t]$. For any $\e\in \{ \pm 1\}^m$, substitute 
$t=(\e_1t_1,\ldots,\e_mt_m)$ in the equation (\ref{P1divQ1}). Then we obtain 
\[ Q_1(t_1^2,t_2^2,\ldots,t_m^2) = P_{\e}(t)g(\e_1t_1,\ldots,\e_mt_m).\]
Thus the polynomial $Q_1(t_1^2,t_2^2,\ldots,t_m^2)$ is divided by $P_{\e}(t)$
 for any $\e\in \{ \pm 1\}^m$, which implies 
$P(t)\mid Q_1(t_1^2,t_2^2,\ldots,t_m^2)$. Therefore 
we conclude that $P(t)$ is irreducible in $\C[t_1^2,t_2^2,\ldots,t_m^2]$. 

Next we prove that the polynomial 
\[ \prod_{\varepsilon_i \in \{+1,-1\}}   \left( 1 + \varepsilon_1 \sqrt{x_1}+\cdots+\varepsilon_m \sqrt{x_m}\right)=\prod_{\e\in \{ \pm 1\}^m}P_\e\left( \sqrt{x_1},\ldots,\sqrt{x_m}\right)
\in \C[x]\]
is irreducible in $\C[x]$. Suppose that 
\[ \prod_{\e\in \{ \pm 1\}^m}P_\e\left( \sqrt{x_1},\ldots,\sqrt{x_m}\right)
=q(x)r(x)\quad (q(x),r(x)\in \C[x]).\]
Substitute $x=(t_1^2,t_2^2,\ldots,t_m^2)$ in the above, then we have 
the equation 
\[ P(t)=q(t_1^2,t_2^2,\ldots,t_m^2)r(t_1^2,t_2^2,\ldots,t_m^2).\]
in $\C[t]$. 
Hence we may assume that $q(t_1^2,t_2^2,\ldots,t_m^2)$ is a constant, so is  
$q(x)\in \C[x]$. Therefore the proof is completed. \QED

\section{The singular locus in the complex torus}
\label{sec:inTorus}

If we can show that the set $\{L_1, \ldots, L_m\}$ is 
a set of  generators of the characteristic ideal
${\rm in}_\www (I(m))$, then we conclude that the sinulgar locus agrees with
the zero set of (\ref{eq:singLocus}).
However, it seems not to be easy to prove it.
Instead of trying to prove it, we will firstly determine the characteristic variety 
in the complex torus.
 
We consider the left ideal $I'(m)$ generated by 
$$ \ell'_i = y_i \theta_i (\theta_i-c_i+1) - (\theta-a)(\theta-b),
\quad i=1, \ldots, m. $$
Here, $\theta_i = y_i \pd{y_i}$ and
$\theta = \theta_1 + \cdots + \theta_m$.
These operators are obtained by
applying the change of the coordinates
$ y_i = 1/x_i $, $i=1, \ldots, n$
to $\ell_i$'s and multiplying $y_i$ to them.
The ring of differential operators with respect to the variable $y$
is also denoted by $D$ as long as no confusion arises.
The characteristic varieties of $I(m)$ and $I'(m)$ agree
on the complex torus $\left( \C^* \right)^m$
under the change of coordinates $y_i = 1/x_i$.

We use the order $\succ_w$ defined by the first weight vector $w^{(1)}=
(\bm{0},\bm{1})$ and the second weight vector $w^{(2)}=(\bm{1},\bm{0})$. In 
other words, $y^\alpha\partial^\beta\succ_w y^{\alpha'}\partial^{\beta'}$ 
if and only if 
\begin{eqnarray*}
&& (\alpha,\beta)\cdot w^{(1)}>(\alpha',\beta')\cdot w^{(1)}\\
&\mathrm{or}& \left( (\alpha,\beta)\cdot w^{(1)}=(\alpha',\beta')\cdot w^{(1)}
\mathrm{~and~} (\alpha,\beta)\cdot w^{(2)}>(\alpha',\beta')\cdot w^{(2)}\right)
\\
&\mathrm{or}& \left(  (\alpha,\beta)\cdot w^{(j)}=(\alpha',\beta')\cdot w^{(j)}
~ (j=1,2) \mathrm{~and~} (\alpha,\beta)>_{\mathrm{lex}}(\alpha',\beta')
\right). 
\end{eqnarray*}
We denote by $\mathrm{in}_{\prec_w}(f)$ the leading monomial of 
$f \in D$
with respect to the order $\prec_w$. 
For two elements $f,g$ $\in D$ with 
\[ \mathrm{in}_{\prec_w}f=f_{\alpha\beta}y^\alpha\xi^\beta,\qquad 
\mathrm{in}_{\prec_w}g=g_{\alpha'\beta'}y^{\alpha'}\xi^{\beta'},
\] 
we define the $S$-pair $\mathrm{sp}(f,g)$ of $f$ and $g$ by
\begin{equation}
\mathrm{sp}(f,g) = g_{\alpha'\beta'}y^{\gamma-\alpha}\partial^{\delta-\beta}f
- f_{\alpha\beta}y^{\gamma-\alpha'}\partial^{\delta-\beta'}g,
\nonumber
\end{equation}
where 
\[
\gamma = \left( \max\{ \alpha_1,\alpha_1'\}, \ldots, \max\{ \alpha_m,
\alpha_m'\}\right),~ 
\delta = \left( \max\{ \beta_1,\beta_1'\}, \ldots, \max\{ \beta_m,\beta_m'\}
\right). 
\]

For a subset $G$ of $D$, the relation 
$f = \sum c_i g_i,  g_i \in G$ is called a standard representation of $f$
with respect to $G$
when $c_i g_i \preceq_{w} f$ holds for all $i$ such that $c_i \not= 0$.

\begin{proposition} \label{prop:charOnTorus}
The characteristic ideal ${\rm in}_\www (I'(m))$
is generated by
${\rm in}_\www(\ell'_i)$, $i=1, \ldots, m$.
\end{proposition}

{\it Proof}\/.
We use the order $\succ_w$. 
Since 
$\mathrm{in}_{\prec_w}(\ell_i')=y_i^3\xi_i^2$, 
we have 
\begin{eqnarray*}
\mathrm{sp}(\ell'_i,\ell'_j) &=& y_j^3\partial_j^2\ell_i'-
y_i^3\partial_i^2\ell_j'. 
\end{eqnarray*}

It is expressed as follows: 
\begin{eqnarray}
&& 
\mathrm{sp}(\ell'_i,\ell'_j)\nonumber\\
 &=& 
(y_j^3\partial_j^2-\ell_j')\ell_i'-(y_i^3\partial_i^2-\ell_i')\ell_j'
-(\ell_i'\ell_j'-\ell_j'\ell_i')\nonumber\\
&=& (y_j^3\partial_j^2-\ell_j')\ell_i'-(y_i^3\partial_i^2-\ell_i')\ell_j'
-(2\theta-a-b+1)(\ell_i'-\ell_j')\nonumber\\
&=& \left\{ y_j^3\partial_j^2-\ell_j'-(2\theta-a-b+1)\right\}\ell_i'\nonumber \\
& & \hspace{3cm} -
\left\{ y_i^3\partial_i^2-\ell_i'-(2\theta-a-b+1)\right\}\ell_j'.
\label{standard_rep}
\end{eqnarray}
Note that we have used the commutation relation
\[\ell_i'\ell_j'-\ell_j'\ell_i'
=-(2\theta-a-b+1)(\ell_i'-\ell_j'), \]
which is obtained by a straightforward calculation. 
By using the relations
\[
\partial_i^2\theta_i = 2\partial_i^2+\theta_i\partial_i^2,\quad
\partial_i^2\theta_i^2 = 
4(1+y_i\partial_i)\partial_i^2+\theta_i\partial_i^2,
\]
we have 
\begin{eqnarray*}
y_j^3\partial_j^2\ell_i' &=& y_i^3y_j^3\partial_i^2\partial_j^2
 -y_j^3\left\{ (-c_i+2)y_i\partial_i+\sum_{k=1}^m\theta_k^2+4(1+y_j\theta_j)
\right.\\
&& \left.
+\sum_{k\neq k'}\theta_k\theta_{k'}+4\sum_{k\neq j}\theta_k
-(a+b)\left( \sum_{k=1}^m \theta_k+2\right)+ab
\right\}\partial_j^2,
\end{eqnarray*}
which implies 
\begin{eqnarray*}
&&\mathrm{in}_{\prec_w}\mathrm{sp}(\ell_i',\ell_j')\nonumber\\
 &=& 
\mathrm{in}_{\prec_w}\left\{ -y_j^3\left( \sum_{k,k'}
\hspace{-0.0cm} y_ky_{k'}\partial_k
\partial_{k'}\right)\partial_j^2
+y_i^3\left( \sum_{k,k'}\hspace{-0.0cm} y_ky_{k'}\partial_k
\partial_{k'}\right)\partial_i^2
\right\}\label{Utsumi}\\
&=& y_1^2y_i^3\xi_1^2\xi_i^2 \label{Sawamura}
\end{eqnarray*}
for $i<j$. 
On the other hand, since
\begin{eqnarray*}
&& y_i^3\partial_i^2-\ell_i'-(2\theta-a-b+1)\\
&=& (c_i-2)y_i\theta_i+(\theta-a)(\theta-b)-(2\theta-a-b+1), 
\end{eqnarray*}
we have 
\begin{eqnarray*} 
&&
 \mathrm{in}_{\prec_w}\left\{ y_i^3\partial_i^2-\ell_i'-(2\theta-a-b+1)\right\}
\\
&=& \mathrm{in}_{\prec_w}\left\{
(c_i-2)y_i\theta_i+(\theta-a)(\theta-b)-(2\theta-a-b+1) \right\}\\
&=& \mathrm{in}_{\prec_w}(\theta^2)\\
&=& y_1^2\xi_1^2. 
\end{eqnarray*}
Note that it is independent of the index $i$. 
Hence we conclude that 
\begin{eqnarray*}
\mathrm{in}_{\prec_w}\left\{ y_j^3\partial_j^2-\ell_j'-(2\theta-a-b+1) 
\right\}\ell_i' &=& y_1^2y_i^3\xi_1^2\xi_i^2,\\
\mathrm{in}_{\prec_w}\left\{ y_i^3\partial_i^2-\ell_i'-(2\theta-a-b+1) 
\right\}\ell_j' &=& y_1^2y_j^3\xi_1^2\xi_j^2, 
\end{eqnarray*}
which imply that
the expression (\ref{standard_rep}) 
is a standard representation of $\mathrm{sp}(\ell'_i,\ell'_j)$
with respect to the set $\{\ell_1'\ldots,\ell_m'\}$ and
the order $\prec_\www$. 
It follows from the Buchberger's criterion that it is a Gr\"obner basis 
with respect to that order.  
Therefore the set 
\[\{ \mathrm{in}_\www(\ell_i')\mid i=1,\ldots,m\}=
\left\{ \left. y_i(y_i\xi_i)^2-\left( \sum_{j=1}^my_j\xi_j\right)^2\hspace{0.1cm}
\right|\hspace{0.1cm} i=1,\ldots,m\right\}\]
is a Gr\"obner basis of $\mathrm{in}_\www(I'(m))$ by the theorem stated in 
\cite[section 2]{oaku-char}
(the condition on the order can be weakened as in 
\cite[Th 1.1.6]{SST}). 
In particular, it is a set of generators of the characteristic ideal 
$\mathrm{in}_\www(I'(m))$. 
\QED

Let us determine the singular locus of $\ell'_1, \ldots, \ell'_m$.
The principal symbol $L_i'$ of $\ell_i'$ is equal to 
\[ L_i'=y_i^3\xi_i^2-\left( \sum_{j=1}^my_j\xi_j\right)^2.\]
When $y_i\neq 0$, it is factored as 
\[ L_i'=\left( y_i\sqrt{y_i}\xi_i-\sum_{j=1}^my_j\xi_j\right)
\left( y_i\sqrt{y_i}\xi_i+\sum_{j=1}^my_j\xi_j\right).\]
We can show that the determinant of the coefficient matrix 
of the system 
\begin{equation}
\nonumber
y_i\sqrt{y_i}\xi_i+\e_i\sum_{j=1}^my_j\xi_j=0,\quad i=1,\ldots,m 
\end{equation}
is equal to 
\begin{equation} \label{eq:singLocus2}
\displaystyle\left(\prod_{j=1}^m y_j\sqrt{y_j}\right)\left( 
 1+\sum_{j=1}^m\dfrac{\e_j}{\sqrt{y_j}}\right)
\end{equation} 
by an analogous way to the proof of the Proposition \ref{prop:det}.
Therefore,
the singular locus of $I'(m)$ is equal to the union of the 
zero sets of (\ref{eq:singLocus2}) 
where $\varepsilon_j$'s run over $\{-1,+1\}$. 
Thus, we have the following theorem.
\begin{theorem}
The singular locus of $I(m)$
agrees with the zero set of (\ref{eq:singLocus}) in the complex torus.
\end{theorem}

\section{Singular locus and the coordinate hyperplanes}

We will prove that the coordinate hyperplanes are contained in 
the singular locus of $I(m)$
by discussing the cohomological solution sheaf
${\cal E}xt^1_{{\cal D}^{an}}({\cal D}^{an}/{\cal D}^{an} I(m), {\cal O}^{an})$.
We need a set of generators of the syzygies of $I(m)$ 
to describe the first cohomological solutions
(as to an algorithmic method to determine it, see, e.g., \cite{takayama-cohom}).
We utilize a Gr\"obner basis with the order $\succ_\km$,
where $\km = (-1, \ldots, -1, 1, \ldots, 1)$
and
the lexicographic order $\succ$
$\pd{1} \succ \cdots \succ \pd{m} \succ x_1 \succ \cdots \succ x_m$
as the tie-breaker,
to determine the syzygies among generators of $I(m)$.

In order to use the $S$-pair criterion,
we will work in the homogenized Weyl algebra 
$D^{(h)}=\C[h]\langle x_1, \ldots, x_m, \pd{1}, \ldots, \pd{m} \rangle$ 
(see, e.g., \cite[\S 1.2]{SST}). 
The variable $h$ is the homogenization variable which commutes with all
other variables and we have the relation $\pd{i} x_i = x_i \pd{i} + h^2$.

Put 
\begin{eqnarray*}
S_i    &=& \theta_i (\theta_i + (c_i-1)h^2), \\
S_{ab} &=& \left(\sum_{i=1}^m \theta_i + a h^2 \right)
           \left(\sum_{i=1}^m \theta_i + b h^2 \right)
\end{eqnarray*}
and
\begin{eqnarray*}
T_i    &=& h S_i - x_i S_{ab}, \\
T_{ij} &=& x_j S_i - x_i S_j.
\end{eqnarray*}
They are homogeneous elements in the $D^{(h)}$.
The operator $T_i$ is the homogenization of $\ell_i$ and
the operator $T_{ij}$ is the homogenization of $x_j \ell_i - x_i \ell_i$.
For two elements in $D^{(h)}$, their $S$-pair with respect to the order 
$\succ_\km$ is defined similarly as in the section \ref{sec:inTorus}. 
We also use the terminology ``standard representation'' analogously 
for elements in $D^{(h)}$. 

\begin{theorem}  \label{th:spair}
The set $G=\{ T_1, \ldots, T_m, T_{12}, T_{13}, \ldots, T_{m-1,m} \}$
satisfies the $S$-pair criterion in the homogenized Weyl algebra $D^{(h)}$;
$G$ is a Gr\"obner basis with respect to the order $\succ_\km$.
\end{theorem}

{\it Proof}\/.
We have the following standard representations of $S$-pairs in terms of $G$: 
\begin{eqnarray}
\mathrm{sp}(T_i,T_j) &=& S_jT_i-S_iT_j ~=~ S_{a-1,b-1}T_{ij},
\label{eq:titj}\\
\mathrm{sp}(T_i,T_{ij}) &=& x_jT_i-hT_{ij} ~=~ x_iT_j,\label{eq:titij}\\
\mathrm{sp}(T_j,T_{ij}) &=& x_i^2\partial_i^2T_j-hx_j\partial_j^2T_{ij}
\nonumber\\
&=& \left\{ x_i(x_j^{-1}S_j)-c_ih^2\theta_i\right\}T_j-(2h^2\theta_j+c_jh^4)
T_i,\\
\mathrm{sp}(T_k,T_{ij}) &=& x_i^2x_j\partial_i^2T_k-hx_k^2\partial_k^2T_{ij}
\nonumber\\
&=& hS_jT_{ki}+x_kS_iT_j-c_ih^2x_j\theta_iT_k+c_kh^3\theta_kT_{ij}, \\
\mathrm{sp}(T_{ij}, T_{ik}) &=& x_kT_{ij}-x_jT_{ik} ~=~ -x_iT_{jk},\\
\mathrm{sp}(T_{ij}, T_{kj}) &=& 
x_k^2\partial_k^2T_{ij}-x_i^2\partial_i^2T_{kj} =
S_jT_{ik}-c_kh^2\theta_kT_{ij}+c_ih^2\theta_iT_{kj},\\
\mathrm{sp}(T_{ij},T_{jk}) &=& 
x_jx_k\partial_j^2T_{ij}-x_i^2\partial_i^2T_{jk}\nonumber\\
&=& \left\{ S_k+(2-c_j)h^2x_k\partial_j\right\} T_{ij}+(c_j-2)h^4T_{ik}
\nonumber\\
&& \hspace{8pt}
+\left\{ (2-c_j)h^2x_i\partial_j+c_ih^2\theta_i-x_i\theta_j\partial_j
\right\} T_{jk}, \\
\mathrm{sp}(T_{ij},T_{i'j'}) &=& 
x_{i'}^2x_{j'}\partial_{i'}^2T_{ij}-x_i^2x_j\partial_i^2T_{i'j'}\nonumber\\
&=& x_{j'}S_jT_{ii'}-x_{i'}S_iT_{jj'}
-c_{i'}h^2x_{j'}\theta_{i'}T_{ij}+c_ih^2x_j\theta_iT_{i'j'}, 
\end{eqnarray}
where we assume that the indices $i,j,k,i',j'$ satisfy $i\neq k$, 
$j\neq k$ and $\{ i,j\}\cap \{ i',j'\}=\phi$. 
Note in the above that we regard $T_{ji}=-T_{ij}$ for $i<j$. 
Thus, we have proved the set $G$ is a Gr\"obner basis. 
\QED


By \cite[Th 9.10]{oaku-takayama}, 
syzygies are generated by the dehomogenizations of the standard representations
of the $S$-pairs. The following Corollary will be used to complete the proof of 
our main theorem.

\begin{corollary}  \label{th:syzygy}
The set of relations
which are derived from the standard representations of the s-pairs
gives a set of generators of the syzygies among
$\ell_i$, $(i=1, \ldots, m)$, $\ell_{ij}=x_j \ell_i - x_i \ell_j$, $1 \leq i < j \leq m$.
For example, {\rm (\ref{eq:titij})} yields the syzygy
$ x_j \ell_i - \ell_{ij} - x_i \ell_j =0$.  
\end{corollary}

\begin{theorem}  \label{th:singLocus}
The singular locus of $I(m)$ is the zero set of
(\ref{eq:singLocus}).
\end{theorem}

{\it Proof}\/.
It follows from the discussions in the section 
2 that
we may prove only that $x_i=0$ are contained in the singular locus.
Let us consider if $x_m=0$ is a singular locus or not.
Let $g_m(x')$ be a non-zero solution of $I(m-1)$.
This function does not depend on $x_m$.
Put $g_1= \cdots = g_{m-1}=0$.
Then, we have $\ell_i \cdot g_m = 0$ for $i\not=m$ and
$\ell_j \cdot g_i=0$ for $i=1, \ldots, m-1$.
Define $g_{ij} = x_j g_i - x_i g_j$.
$\sum g_i e_i + \sum g_{ij} e_{ij}$ are annihilated by the generators of the syzygies
given in the Corollary \ref{th:syzygy}. 
For instance, we have the syzygy 
\[(\theta_j(\theta_j-1)+c_j\theta_j)\ell_i-(\theta_i(\theta_i-1)+c_i\theta_i)
\ell_j - (\theta+a-1)(\theta+b-1)\ell_{ij}=0\]
by the equation (\ref{eq:titj}). For $i<j=m$, we have 
\begin{eqnarray*}
&& (\theta_j(\theta_j-1)+c_j\theta_j)g_i-(\theta_i(\theta_i-1)+c_i\theta_i)
g_j - (\theta+a-1)(\theta+b-1)g_{ij}\\
&=& -(\theta_i(\theta_i-1)+c_i\theta_i)g_m-(\theta+a-1)(\theta+b-1)(-x_i g_m)\\
&=& -(\theta_i(\theta_i-1)+c_i\theta_i)g_m+x_m(\theta+a)(\theta+b)g_m\\
&=& -\left\{ 
(\theta_i(\theta_i-1)+c_i\theta_i)-x_m\left(\sum_{k=1}^{m-1}\theta_k
+a\right)
\left(\sum_{k=1}^{m-1}\theta_k
+b\right)\right\}g_m\\
&& \hspace{4pt}(\mathrm{since~} g_m=g_m(x')~ \mathrm{does~ not~ depend~ on~} 
x_m, \mathrm{we~ have~} \partial_mg_m=0)\\
&=& 0. 
\end{eqnarray*}
When $i<j<m$, obviously the equation 
\[(\theta_j(\theta_j-1)+c_j\theta_j)g_i-(\theta_i(\theta_i-1)+c_i\theta_i)
g_j - (\theta+a-1)(\theta+b-1)g_{ij}=0\]
holds. 

Let us try to solve
$\ell_i \cdot f = g_i$, $i=1, \ldots, m$
and 
$\ell_{ij} \cdot f = g_{ij}$, $1 \leq i < j \leq m$.
The second group of equation is solved when the first group is solved.
Put $f=\sum_{k=0}^\infty f_k x_m^k$.
The left hand side of $\ell_m f$ can be factored by $x_m$.
On the other hand, the right hand side $g_m$ is nonzero and does not
depend on $x_m$.
Therefore the system $\ell_i \cdot f=g_i$ does not have a holomorphic solution
along $x_m=0$.
Therefore, we have proved that 
${\cal E}xt^1_{\cal D}({\cal D}^{an}/{\cal D}^{an}E_C, {\cal O})$ 
is not zero at a generic point in $x_m=0$.
By Kashiwara's theorem 
[8,Theorem 4.1],
the ${\cal E}xt^1$ must be zero if $x_m=0$ is not a singular locus.
Thus, we have proved that $x_m=0$ is the singular locus.
We can analogously show that other varieties $x_i=0$ are also 
contained in the singular locus.
\QED

\section{The $A$-hypergeometric system associated to the Lauricella $F_C$}

The binomial $D$-modules \cite{binomial-d-modules}
are introduced to study
classical hypergeometric systems including the Lauricella $F_C$.
The contents of the first half part of this section are implicitly or explicitly explained in
\cite{binomial-d-modules},
but they do not seem to be publicized to people
who study classical Lauricella functions and related topics.
We add the first part of this section to explain how to apply the theory of $A$-hypergeometric
systems and binomial $D$-modules to study $F_C$.
The second part contains a new result and utilizes the first part;
The last Theorem \ref{th:a-sing} describes
the singular locus for the $A$-hypergeometric system
associated to $F_C$ in the complex torus.
The singular locus is the zero set of the principal $A$-determinant \cite{GKZ-det}
for the $A$ associated to $F_C$.

We denote $I(m)$, which annihilates the Lauricella function $F_C$,
by $E_C$ in this section.
For a given Heun system, there exists 
a corresponding binomial $D$-module. 
In case of $E_C$, the corresponding binomial system is 
an $A$-hypergeometric system.
Let us study this system.

Let $e_1, \ldots, e_{m+1},e_{m+2}$ be the standard basis of
$\Z^{m+2}$.
Following \cite{S1995},
consider the set of points ${\cal A}=\{ e_1 + e_{m+2}, e_2+e_{m+2}, \ldots, e_{m+1}+e_{m+2},
-e_1 + e_{m+2}, -e_2+e_{m+2}, \ldots, -e_{m+1}+e_{m+2}\}$.
We define a matrix $A(F_C,m)$ consisting of these points as column vectors.
The matrix is an $(m+2) \times 2(m+1)$ matrix.
For example, we have
$$
A(F_C,2)=
\left(
\begin{array}{cccccc}
 1& 0& 0&  - 1& 0& 0 \\
0&  1& 0& 0&  - 1& 0 \\
0& 0&  1& 0& 0&  - 1 \\
 1&  1&  1&  1&  1&  1 \\
\end{array}
\right)
$$
Let $H_A(\beta)$ be the $A$-hypergeometric system
associated to the matrix $A(F_C,m)$,
the parameter
$\beta^T=(1-c_1, \ldots, 1-c_m, b-a, \sum_{j=1}^m c_j-a-b-m)$
and the variables 
$u_1, \ldots, u_{m+1},
u_{-1}, \ldots, u_{-(m+1)}$
as independent variables.
The associated differential operators for $u_j$ and $u_{-j}$
are denoted by $\pd{j}$ and $\pd{-j}$ respectively.
The toric ideal $I_A$ defined by $A=A(F_C,m)$,
$I_A = \{ \pd{}^u-\pd{}^v \,|\, Au=Av, u,v \in \N_0^{2m+2} \}$,
is generated by
$\pd{j} \pd{-j} - \pd{m+1} \pd{-(m+1)}$,
$j=1, \ldots, m$
in $\C[\pd{1}, \ldots, \pd{m+1}, \pd{-1}, \ldots, \pd{-m-1}]$.

The left ideal $H_A(\beta)$ is generated by the row vectors
of $A \theta_u -\beta$
and $I_A$
where $\theta_u = (u_1 \pd{1}, \ldots, u_{m+1}\pd{m+1}, 
                  u_{-1}\pd{-1}, \ldots, u_{-(m+1)}\pd{-(m+1)})^T$.
We denote the $i$-th row vector of $A \theta - \beta$ by $E_i-\beta_i$.

\begin{theorem} {\rm (\cite{binomial-d-modules})}
${\rm Sol}(E_C) \simeq {\rm Sol}(H_A(\beta))$
holds for any parameter
in a simply connected neighborhood of any generic point.
In particular, the holonomic rank of $E_C=I(m)$ is equal to $2^m$.
\end{theorem}

{\it Proof}\/.
Let $F$ be a solution of $E_C$.
Following \cite{S1995}, we consider the following function
\begin{equation}  \label{eq:FC2A}
f(u)=
u_{m+1}^{-a} u_{-(m+1)}^{-b} \prod_{j=1}^m u_{-j}^{c_j-1}
 F \left( 
  \frac{u_{1} u_{-1}}{u_{m+1} u_{-(m+1)}}, \ldots,
  \frac{u_{m} u_{-m}}{u_{m+1} u_{-(m+1)}}
 \right)
\end{equation}
Let us prove that the function $f(u)$ is a solution of $H_A(\beta)$.
It is easy to see that $(E_i-\beta_i) \cdot f = 0$.
Put $\eta=u_{m+1}^{-a} u_{-(m+1)}^{-b} \prod_{j=1}^m u_{-j}^{c_j-1}$
and
$z_j =   \frac{u_{j} u_{-j}}{u_{m+1} u_{-(m+1)}}$.
We have
\begin{eqnarray}
\theta_j \theta_{-j} \cdot f(u)
&=& \theta_j \cdot \left( (c_j-1) \eta F + \eta z_j F_j \right) \nonumber \\
&=& (c_j-1) \eta z_j F_j + \eta z_j F_j+\eta z_j z_j F_{jj} \nonumber \\
&=& \eta \, \theta_{z_j} ( \theta_{z_j} + c_j-1) \cdot F(z) \label{eq:FC2A1}
\end{eqnarray}
where $F_j$ denotes the partial derivative of $F(z_1, \ldots, z_m)$
with respect to the variable $z_j$.
Analogously, we get
\begin{equation} \label{eq:FC2A2}
\theta_{m+1} \theta_{-(m+1)} \cdot f(u)
= \eta
  \left( \sum_{i=1}^m \theta_{z_i} + a \right)
  \left( \sum_{i=1}^m \theta_{z_i} + b \right) \cdot F(z)
\end{equation}
Put 
${\tilde \ell}_j = 
u_ju_{-j}u_{m+1} u_{-(m+1)} ( \pd{j} \pd{-j} - \pd{m+1} \pd{-(m+1)} )$.
This is equal to 
$$ 
  u_{m+1} u_{-(m+1)} \theta_j \theta_{-j}-
  u_ju_{-j} \theta_{m+1} \theta_{-(m+1)}.
$$
It follows from (\ref{eq:FC2A1}) and (\ref{eq:FC2A2})
that we have ${\tilde \ell}_j \cdot f(u) = 0$,
which implies that $f(u)$ is a solution of $H_A(\beta)$.
Note that the correspondence from $F$ to $f$ is an injection
among $\C$-vector spaces
on a simply connected open set.

Conversely, let $f$ be a solution of $H_A(\beta)$.
We define new $2m+2$ variables $z_j$ by
\begin{eqnarray}  \label{eq:birational}
z_j &=& u_j u_{-j}/(u_{m+1} u_{-(m+1)}), \ 
z_{m+j} = u_{-j}^{-1}, \ j=1, \ldots, m, \\
z_{2m+1} &=& u_{m+1}^{-1},
z_{2m+2} = u_{-(m+1)}^{-1}  \nonumber
\end{eqnarray}
Note that this gives an isomorphism of the complex toruses
$(\C^*)^{2m+1} = \{ u \}$ and
$(\C^*)^{2m+1} = \{ z \}$.
The Euler operator $\theta_{\pm u_j}=u_{\pm j} \pd{\pm j}$ 
can be written as a sum of Euler operators with respect to $z_i$'s.
In fact, we have
\begin{eqnarray*}
 u_j \pd{j} &=& z_j \pd{z_j}, u_{-j} \pd{-j} = z_j \pd{z_j} - z_{m+j} \pd{z_{m+j}} \\
 u_{m+1} \pd{m+1} &=& -\sum_{k=1}^m z_j \pd{z_k} - z_{2m+1} \pd{z_{2m+1}} \\
 u_{-(m+1)} \pd{-(m+1)} &=& -\sum_{k=1}^m z_j \pd{z_k} - z_{2m+2} \pd{z_{2m+2}}.
\end{eqnarray*}
Put $f'=\eta^{-1} f$.
The equations $(E_i -\beta_i) \cdot f = 0$ yield
$ \theta_{z_j} \cdot f' = 0$ for $j=m+1, \ldots, 2(m+1)$.
This implies that $f'$ depends only on $z_1, \ldots, z_m$.
An analogous calculation with (\ref{eq:FC2A1}) and (\ref{eq:FC2A2})
yields the equation $\ell_i \cdot f'(u(z)) = 0$. 
This means that the map $F(z) \mapsto f(u)$ is surjective.
Thus, we have proved
${\rm Sol}(E_C) \simeq {\rm Sol}(H_A(\beta))$.

 
The correspondence gives the holonomic rank of $E_C$
by  evaluating the degree of $I_A$ \cite{GKZ}.
Since $\{ \underline{\pd{j} \pd{-j} } - \pd{m+1} \pd{-(m+1)} \,|\,
 j=1, \ldots, m \}$ is a Gr\"obner basis of $I_A$,
the degree is equal to that of the monomial ideal
generated by $\pd{j} \pd{-j}$, $j=1, \ldots, m$.
This degree is equal to $2^m$.
\QED

An interesting application of this isomorphism is
the following irreducibility condition of $E_C$.
We can utilize the recent result by Beukers \cite{B2010}  and
Schulze and Walther \cite{SW2010}
on irreducibility of $A$-hypergeometric systems
to give a condition of the irreducibility of $E_C$.
\begin{theorem} {\rm (\cite{B2010}, \cite{SW2010})}
The system $E_C$ is irreducible if and only if
$$\frac{1}{2}\left(
  \sum_{i=1}^m c_i - a -b -2 
  \sum_{i=1}^m \varepsilon_i (1-c_i) + \varepsilon_{m+1} (b-a) \right)
 \not\in \Z
$$
for all combinations of $\varepsilon_i \in \{ -1, 1 \}$.
\end{theorem}

{\it Proof}\/.
It follows from the previous theorem that the irreducibility of $E_C$
is equivalent to that of $H_A(\beta)$.
In fact, the solution space is locally isomorphic and
differential operators with rational function coefficients in $z$
is mapped to those in $u$ by (\ref{eq:birational}).
The primitive integral support functions $P_J(s)$ for $A(F_C,m)$ are
$(1/2)(s_{m+2}+\sum_{j \in J} s_j-\sum_{j \not\in J} s_j)$,
$ J \subseteq [1,m+1]$ 
where $s_j$'s are the dual basis for the $e_i$'s
\cite{S1995}.
It follows from \cite{B2010} or \cite{SW2010} that
the irreducibility condition is that $P_J(\beta) \not\in \Z$ for all $J$,
which is equivalent to the condition in the theorem.
\QED


Finally, we discuss the singular locus of the $H_A(\beta)$
via the correspondence.
The correspondence is not only for the classical solutions as we have seen, 
but also for some
$D$-module invariants including the singular locus on the complex torus.
In this case, we utilize our result on the singular locus for $F_C$ 
to derive a result on the $A$-hypergeometric system.
\begin{theorem} \label{th:a-sing}
The singular locus of $H_A(\beta)$ in the complex torus is given
by  the zero set of 
$$
  \prod_{\varepsilon_i \in \{ -1, 1 \} } 
    \left( 1 + \varepsilon_1 \sqrt{\frac{u_1 u_{-1}}{u_{m+1}u_{-(m+1)}}}
             +\cdots
             +\varepsilon_m \sqrt{\frac{u_m u_{-m}}{u_{m+1}u_{-(m+1)}}}
    \right)
$$
\end{theorem}

{\it Proof}\/.
We denote by $D^*_{2m+2}$ the ring of differential operators on the complex torus
$$ \C\langle z_1^{\pm}, \ldots, z_{2m+2}^{\pm}, \pd{z_1}, \ldots, \pd{z_{2m+2}} \rangle.$$
Let $I$ be a left ideal in $D^*_{2m+2}$ and for a complex number $\alpha$
denote by
$D^* z_{2m+2}^\alpha$
the left $\C\langle z_{2m+2}^{\pm} \pd{z_{2m+2}}\rangle$-module
$\C\langle z_{2m+2}^{\pm} \pd{z_{2m+2}}\rangle/\langle z_{2m+2} \pd{z_{2m+2}} - \alpha\rangle$.
The outer tensor product 
$(D^*_{2m+2}/I) \outertensor D^* z_{2m+2}^\alpha$ 
is defined by the restriction of
$$D^*_{2m+3}/\langle I, z_{2m+3} \pd{z_{2m+3}}-\alpha \rangle$$
to $z_{2m+3}-z_{2m+2}=0$.
In other words,
\begin{eqnarray} \label{eq:outertensor}
& &(D^*_{2m+2}/I)\outertensor D^* z_{2m+2}^\alpha
 |_{z_{2m+2} \mapsto t, \pd{z_{2m+2}} \mapsto \pd{t} }  \\
&\simeq&
D^*_{2m+1}\langle t^{\pm},s,\pd{t},\pd{s}\rangle/
\left( 
 \langle I,-t\pd{s}+s\pd{s}-\alpha \rangle 
 + sD^*_{2m+1}\langle t^{\pm},s,\pd{t},\pd{s}\rangle
\right)  \nonumber
\end{eqnarray}
where we make replacements
$$ s = z_{2m+2}-z_{2m+3}, t = z_{2m+2},
 z_{2m+2} \pd{2m+2} = t \pd{s} + t \pd{t},
 z_{2m+3} \pd{2m+3} = -t \pd{s} + s \pd{s}
$$
in $I$.
Let $u=(0, \ldots, 0,1)$ be the weight vector where
$1$ stands for the variable $s$.
Then,
$b={\rm in}_{(-u,u)}(-t\pd{s}+s \pd{s}-\alpha) = -t \pd{s}$.
Since $t$ is invertible, we have
$ D^*_{2m+1}\langle t^{\pm},s,\pd{t},\pd{s} \rangle b \cap \C[s \pd{s}]
= \langle s \pd{s} \rangle$.
Therefore, by the restriction algorithm (see, e.g., \cite{oaku-takayama}),
we can prove that
(\ref{eq:outertensor}) is isomorphic to
$$ D^*_{2m+1} \langle t^{\pm} ,\pd{t} \rangle/
\left( 
\left(
 \langle I,-t\pd{s}+s\pd{s}-\alpha \rangle 
 + sD^*_{2m+1}\langle t^{\pm},s,\pd{t},\pd{s}\rangle
\right) \cap D^*_{2m+1}\langle t^{\pm}, \pd{t} \rangle
\right).
$$
of which denominator is called the restriction ideal.

Put $\tau_j = z_j \pd{z_j}$.
Let $I$ be the hypergeometric ideal $H_A(\beta)$ expressed in terms of the variable
in $z_j$ (\ref{eq:birational}), 
which is generated in $D^*_{2m+2}$ by
$$
\tau_j (\tau_j - \tau_{m+j}) - z_j
(\sum_{j=1}^m \tau_j + \tau_{2m+1}) (\sum_{j=1}^m \tau_j + \tau_{2m+2}), \quad
j=1, \ldots, m
$$
and
$$ 
\tau_{m+j}-(1-c_j), 
\quad j=1, \ldots, m, \quad
\tau_{2m+2}-\tau_{2m+1}-(b-a),
\tau_{2m+1}-a.
$$
Note that
$\tau_j (\tau_j - \tau_{m+j}) - 
z_j (\sum_{j=1}^m \tau_j + \tau_{2m+1}) (\sum_{j=1}^m \tau_j + t\pd{s} + t \pd{t})$
is in $I$ under the change of variables from $z_{2m+2}, z_{2m+3}$ to $s,t$.
Subtracting $(\sum \tau_j + \tau_{2m+1})(-t\pd{s}+s\pd{s}-\alpha)$ from it,
we conclude that the restriction ideal contains
\begin{equation} \label{eq:restrictionElem}
\tau_j (\tau_j - \tau_{m+j}) - 
z_j (\sum \tau_j + \tau_{2m+1})(\sum \tau_j + t \pd{t}-\alpha).
\end{equation}
We have defined the outer tensor product by $D^* z_{2m+2}^\alpha$ and
studied its properties.
We can make analogous discussions for outer tensor products by other variables
and we conclude from (\ref{eq:restrictionElem}) that
there exists a left ideal $I'$ such that
\begin{equation} \label{eq:outertensorAll}
 (D^*_{2m+2}/I) \outertensor D^* \eta^{-1}
\simeq D^*_{2m+2}/I'
\end{equation}
and 
$I' \supseteq \langle I(m), \pd{z_{m+1}}, \ldots, \pd{z_{2m+2}} \rangle$.
Here, $I(m)$ is the left ideal generated by the Lauricella operators
(\ref{eq:LauricellaOp})
($x_i$'s are replaced by $z_i$'s respectively).

We denote by ${\rm Sing}^*(M)$ the singular locus of $M$ in the complex torus.
It follows from (\ref{eq:outertensorAll}) that
\begin{eqnarray}
& & {\rm Sing}^*(D^*_{2m+2}/I) =
{\rm Sing}^*((D^*_{2m+2}/I) \outertensor D^* \eta^{-1}) \nonumber \\
&=& {\rm Sing}^*(D^*_{2m+2}/I') 
\subseteq {\rm Sing}^*(D^*/I(m))
\end{eqnarray}
Since ${\rm Sing}^*(D^*/I(m))$ is irreducible,
the singular locus of the $A$-hypergeometric system
${\rm Sing}^*(D^*_{2m+2}/I)$ is empty or agrees with ${\rm Sing}^*(D^*/I(m))$.
Since the toric ideal $I_A$ is Cohen-Macaulay when $A=A(F_C,m)$,
the singular locus of the $A$-hypergeometric system $H_A(\beta)$ does not depend
on the parameter $\beta$ \cite{GKZ}, \cite[Section 4.3]{SST}.
Then, we may suppose that $H_A(\beta)$ is irreducible.
By \cite{Hotta}, the $A$-hypergeometric system is regular holonomic and
the irreducibility implies that the irreducibility of the monodromy representation.
If the singular locus in the complex torus is empty, then the monodromy representation
is reducible.
Then, the singular locus is not empty and then we obtain the conclusion.
\QED

\end{document}